\documentclass[11pt,a4paper]{amsart} 
\usepackage[all]{xy}
\SelectTips{cm}{}
\calclayout \makeatletter
\def\serieslogo@{} 
\def\@setcopyright{} 
\makeatother

\title[Thick subcategories of modules over commutative rings]{Thick
subcategories of modules over commutative rings\\ {\footnotesize (with an appendix by
Srikanth Iyengar)}}

\thanks{Version from March 4, 2007.}
\author{Henning Krause}
\address{Srikanth Iyengar\\ Department of Mathematics\\
University of Nebraska\\ Lincoln NE 68588\\ U.S.A.}
\email{iyengar@math.unl.edu}
\address{Henning Krause\\ Institut f\"ur Mathematik\\
Universit\"at Paderborn\\ 33095 Paderborn\\ Germany.}
\email{hkrause@math.upb.de}

\newtheorem{lem}{Lemma}[section]
\newtheorem{prop}[lem]{Proposition}
\newtheorem{cor}[lem]{Corollary}
\newtheorem{thm}[lem]{Theorem}

\theoremstyle{remark}
\newtheorem{rem}[lem]{Remark}

\theoremstyle{definition}
\newtheorem{exm}[lem]{Example}
\newtheorem{defn}[lem]{Definition}

\numberwithin{equation}{section}

\renewcommand{\mod}{\operatorname{mod}\nolimits}

\newcommand{\depth}{\operatorname{depth}\nolimits}

\newcommand{\Ass}{\operatorname{Ass}\nolimits}
\newcommand{\Supp}{\operatorname{Supp}\nolimits}

\newcommand{\Inj}{\operatorname{Inj}\nolimits}

\newcommand{\Mod}{\operatorname{Mod}\nolimits}

\newcommand{\Hom}{\operatorname{Hom}\nolimits}

\newcommand{\Coker}{\operatorname{Coker}\nolimits}

\renewcommand{\dim}{\operatorname{dim}\nolimits}

\newcommand{\Ext}{\operatorname{Ext}\nolimits}
\newcommand{\Tor}{\operatorname{Tor}\nolimits}

\newcommand{\Spec}{\operatorname{Spec}\nolimits}

\newcommand{\height}{\mathrm{ht}}

\newcommand{\lotimes}{\otimes^{\mathbf L}}

\def\C{{\mathcal C}}
\def\D{{\mathcal D}}

\def\bfD{{\mathbf D}}

\def\bbZ{{\mathbb Z}}

\def\fra{{\mathfrak a}}
\def\frm{{\mathfrak m}}
\def\frp{{\mathfrak p}}
\def\frq{{\mathfrak q}}

\begin{document}

\begin{abstract}
For a commutative noetherian ring $A$, we compare the support of a
complex of $A$-modules with the support of its cohomology. This leads
to a classification of all full subcategories of $A$-modules which are
thick (that is, closed under taking kernels, cokernels, and
extensions) and closed under taking direct sums.
\end{abstract}

\maketitle

\section{Introduction}
Let $A$ be a commutative noetherian ring. We consider the category
$\Mod A$ of $A$-modules and the spectrum $\Spec A$ of prime ideals of
$A$. Given a complex $X$ of $A$-modules, we wish to relate its support
(in the sense of Foxby \cite{F}) $$\Supp X=\{\frp\in\Spec A\mid
X\lotimes_A k(\frp)\neq 0\}$$ to the support of its cohomology
$$\Supp H^*X=\bigcup_{i\in\bbZ}\Supp H^iX.$$ In some cases we have the
equality $$\Supp X=\Supp H^*X,$$ for example when $A$ is a Dedekind
domain or when $H^*X$ is finitely generated over $A$. The failure of
this equality is the main theme of recent joint work with Benson and
Iyengar \cite{BIK}, which is motivated by the study of support
varieties of modular representations.  In this paper we establish a
closely related classification of thick subcategories of $\Mod A$ and
prove the following result.

\begin{thm}\label{th:main}
For a subset $\Phi$ of $\Spec A$ the following conditions are equivalent:
\begin{enumerate}
\item For every complex $X$ of $A$-modules we have
$$\Supp X\subseteq \Phi\quad \Longleftrightarrow \quad \Supp H^*X\subseteq \Phi.$$
\item The $A$-modules $M$ with $\Supp M\subseteq \Phi$ form a thick
subcategory of $\Mod A$.
\item Every map $I^0\to I^1$ between injective $A$-modules with $\Ass
I^i\subseteq \Phi$ ($i=0,1$) can be completed to an exact sequence
$I^0\to I^1\to I^2$ such that $I^2$ is injective and $\Ass
I^2\subseteq \Phi$.
\end{enumerate}
\end{thm}

Recall that a classical result of Gabriel \cite{G} provides a
bijection between the set of localizing subcategories of $\Mod A$ and
the set of specialization closed subsets of $\Spec A$. More recently,
a number of authors studied subcategories of $\Mod A$ in terms of
subsets of $\Spec A$; see \cite{GP,H,T}. In this paper we generalize
Gabriel's result in two directions. The first direction (via
associated primes) is fairly elementary but seems to be new. The
second direction (via support) leads to a classification of thick
subcategories of $\Mod A$ and corrects a result of Hovey in \cite{H}.
This classification is formulated in terms of subsets
$\Phi\subseteq\Spec A$ satisfying the equivalent conditions of
Theorem~\ref{th:main}. We call such subsets \emph{coherent} and
establish the following.

\begin{thm}
For a commutative noetherian ring $A$ the following conditions are
equivalent:
\begin{enumerate}
\item The Krull dimension of $A$ is at most one.
\item Every  subset of $\Spec A$ is coherent.
\item $\Supp X=\Supp H^*X$ for every complex $X$ of $A$-modules.
\end{enumerate}
\end{thm}

The proof of this theorem is illustrated by some explict examples of
subsets of $\Spec A$ which are not coherent. It would be interesting
to have a geometric interpretation of coherent subsets in terms of the
Zariski topology on $\Spec A$.

\section{Subcategories via associated primes}

Let $M$ be an $A$-module. Recall that a prime ideal $\frp$ is
\emph{associated} to $M$, if $A/\frp$ is isomorphic to a submodule of
$M$. We denote by $\Ass M$ the set of all prime ideals which are
associated to $M$.

\begin{thm}\label{th:bij}
The map sending a subcategory $\C$ of $\Mod A$ to
$$\Ass\C=\bigcup_{M\in\C}\Ass M$$ induces a bijection between the set
of full subcategories of $\Mod A$, which are closed under taking
submodules, extensions, and direct unions, and the set of subsets of
$\Spec A$. The inverse map sends a subset $\Phi$ of $\Spec A$ to
$$\Ass^{-1}\Phi=\{M\in\Mod A\mid\Ass M\subseteq \Phi\}.$$
\end{thm}

The proof uses some basic facts about associated primes and the
structure of injective modules. The injective envelope of a module $M$
is denoted by $E(M)$ and we observe that $\Ass E(M)=\Ass M$.

\begin{lem}\label{le:ass}
Let $M$ be an $A$-module. Given a submodule $N\subseteq M$ and a
family $(M_i)$ of submodules satisfying $M=\bigcup_i M_i$, we have
$$\Ass N\subseteq\Ass M\subseteq\Ass N\cup \Ass M/N\quad\text{and}\quad
\Ass M=\bigcup_i \Ass M_i.$$
\end{lem}
\begin{proof} 
See \cite[Chap.~IV, \S 1]{B0}.
\end{proof}

For a prime ideal $\frp$ we denote by $k(\frp)$ its residue field.

\begin{lem}\label{le:inj}
Let $I$ be an indecomposable injective $A$-module and $\frp$ its
associated prime ideal. Then $I$ is obtained from  $A/\frp$ by taking
extensions and direct unions. More precisely,
\begin{enumerate} 
\item $I$ is obtained from  $k(\frp)$ by taking
extensions and direct unions, and
\item $k(\frp)$ is a direct union of copies of $A/\frp$.
\end{enumerate}
\end{lem}
\begin{proof}
We sketch the argument and refer to \cite[Chap.~X, \S 8]{B} for
details. For each integer $n\geq 0$ let $I_n$ denote the submodule of
$I$ consisting of all elements annihilated by $\frp^n$. Then we have
$I=\bigcup_{n\geq 0}I_n$ and each factor $I_{n+1}/I_n$ is isomorphic
to a finite direct sum of copies of the residue field $k(\frp)$. Now
observe that $k(\frp)$ is the field of fractions of $A/\frp$ and
therefore a direct union of the form
\begin{equation*}
k(\frp)=\bigcup_{0\neq x\in A/\frp}x^{-1}A/\frp.\qedhere
\end{equation*}
\end{proof}

\begin{proof}[Proof of Theorem~\ref{th:bij}]
Let $\Phi$ be  a subset of $\Spec A$. Then the subcategory $\Ass^{-1}\Phi$
is closed under taking submodules, extensions, and direct unions, by
Lemma~\ref{le:ass}.  Clearly, we have
$$\Ass(\Ass^{-1}\Phi)=\Phi.$$ Now let $\C$ be a subcategory of $\Mod A$, which
is closed under taking submodules, extensions, and direct unions.
We claim that
$$\Ass^{-1}(\Ass\C)=\C.$$ The inclusion $\Ass^{-1}(\Ass\C)\supseteq\C$
is clear. Now suppose that $M$ is a module contained in
$\Ass^{-1}(\Ass\C)$.  Then its injective envelope $E(M)$ is a direct
sum of copies of the form $E(A/\frp)$ with $\frp\in \Ass\C$, since
$\Ass E(M)=\Ass M$. But $\frp\in\Ass\C$ implies $A/\frp\in\C$, and
therefore $E(A/\frp)$ belongs to $\C$, by Lemma~\ref{le:inj}. It
follows that $E(M)$ belongs to $\C$ and therefore $M\in\C$. This
finishes the proof.
\end{proof}

We state a number of consequences of Theorem~\ref{th:bij}.

\begin{cor}
For a full subcategory $\C$ of $\Mod A$ the following conditions are
equivalent.
\begin{enumerate}
\item $\C$ is closed under taking submodules,
extensions, and direct unions.
\item There exists a subset $\Phi$ of $\Spec A$  such that $\C$ consists
of all $A$-modules $M$ satisfying $\Ass M\subseteq \Phi$.
\item There exists an injective $A$-module $I$ such that $\C$ consists
of all $A$-modules which admit a monomorphism into a direct sum of copies of $I$.
\end{enumerate}
\end{cor}
\begin{proof}
(1) $\Leftrightarrow$ (2): This is an immediate consequence of Theorem~\ref{th:bij}.

(2) $\Rightarrow$ (3): Take $I=\bigoplus_{\frp\in \Phi}E(A/\frp)$. Then
$\Ass M\subseteq \Phi$ for every submodule $M$ of a direct sum of copies
of $I$. On the other hand, if $\Ass M\subseteq \Phi$, then $\Ass
E(M)\subseteq \Phi$ and therefore $E(M)$ is a submodule of a direct sum
of copies of $I$.

(3) $\Rightarrow$ (1): Clear.
\end{proof}

Next we restrict the map $\C\mapsto\Ass\C$ to the category $\mod A$ of
all finitely generated $A$-modules.

\begin{cor}[{\cite[Theorem~4.1]{T}}]
The map $\D\mapsto\Ass\D$ induces a bijection between the set
of full subcategories of $\mod A$, which are closed under taking
submodules and extensions, and the set of subsets of $\Spec A$.
\end{cor}
\begin{proof} 
Consider the map $\C\mapsto \C\cap\mod A$ between
\begin{enumerate}
\item[(i)] the set of full subcategories of $\Mod A$, which are closed
under taking submodules, extensions, and direct unions, and
\item[(ii)] the set of full subcategories of $\mod A$, which are
closed under taking submodules and extensions.
\end{enumerate}
This map is bijective; its inverse sends a subcategory $\D$ from (ii)
to the full subcategory of $\Mod A$ consisting of all direct unions of
modules in $\D$. The composition of the first map $\C\mapsto\C\cap\mod
A$ with $\D\mapsto\Ass\D$ is the bjection from
Theorem~\ref{th:bij}. Thus $\D\mapsto\Ass\D$ is bijective.
\end{proof} 

Recall that a full subcategory $\C$ of $\Mod A$ is {\em localizing} if
$\C$ is closed under taking submodules, factor modules, extensions,
and direct sums.  A subset $\Phi$ of $\Spec A$ is {\em specialization
closed} if for any pair $\frp\subseteq\frq$ of prime ideals, $\frp\in
\Phi$ implies $\frq\in \Phi$.

\begin{cor}[{\cite[p.~ 425]{G}}]
The map $\C\mapsto\Ass\C$ induces a bijection between the set
of localizing subcategories of $\Mod A$ and the set of specialization
closed subsets of $\Spec A$.
\end{cor}
\begin{proof}
Suppose that $\C$ is localizing and let $\frp\subseteq\frq$ be prime
ideals with $\frp$ in $\Phi=\Ass\C$. Then $A/\frp\in\C$ and therefore
$A/\frq\in\C$, because $A/\frq$ is a factor module of $A/\frp$. Thus
$\frq\in \Phi$, and we have that $\Phi$ is specialization closed.

Now suppose that $\Phi\subseteq \Spec A$ is specialization closed and let
$N\subseteq M$ be $A$-modules with $M$ in $\C=\Ass^{-1}\Phi$. Then $M/N$
belongs to $\C$, since
$$\Ass M/N\subseteq\{\frp\in\Spec A\mid
(M/N)_\frp\neq0\}\subseteq\{\frp\in\Spec A\mid M_\frp\neq 0\}\subseteq
\Phi$$ where the last inclusion uses that $\Phi$ is specialization
closed. Thus $\C$ is localizing.
\end{proof}

\section{Subcategories via support}
\label{se:supp}
A full subcategory $\C$ of $\Mod A$ is called \emph{thick} if for each
exact sequence $$M_1\to M_2\to M_3\to M_4\to M_5$$ of $A$-modules with
$M_i$ in $\C$ for $i=1,2,4,5$, the module $M_3$ belongs to $\C$.  Note
that a thick subcategory is an abelian category and that the inclusion
functor is exact.

We wish to classify all thick subcategories of $\Mod A$ which are
closed under taking direct sums. For this a few definitions are
needed.

Let $M$ be an $A$-module. Following \cite{F}, the \emph{support} of
$M$ is by definition
$$\Supp M=\{\frp\in\Spec A\mid \Tor_*^A(M,k(\frp))\neq 0\}.$$ For example,
$\Supp I=\Ass I$ for every injective $A$-module $I$.

Let $\Phi$ be a subset of $\Spec A$ and define the full subcategory
$$\Inj_\Phi A=\{I\in\Mod A\mid I \text{ is injective and }\Ass
I\subseteq \Phi\}.$$ We call $\Phi$ \emph{coherent}\footnote{The term
\emph{coherent} refers to the characterizing property of a coherent
ring that every morphism $P_1\to P_0$ between finitely generated
projective modules can be completed to an exact sequence $P_2\to
P_1\to P_0$ such that $P_2$ is finitely generated projective.} if each
morphism $I^0\to I^1$ in $\Inj_\Phi A$ can be completed to an exact
sequence $I^0\to I^1\to I^2$ with $I^2$ in $\Inj_\Phi A$. For example,
each specialization closed subset of $\Spec A$ is coherent.

\begin{thm}\label{th:thick}
The map sending a subcategory $\C$ of $\Mod A$ to
$$\Supp\C=\bigcup_{M\in\C}\Supp M$$ induces a bijection between the
set of full subcategories of $\Mod A$, which are thick and closed
under taking direct sums, and the set of coherent subsets of $\Spec
A$.  The inverse map sends a subset $\Phi$ of $\Spec A$ to
$$\Supp^{-1}\Phi=\{M\in\Mod A\mid\Supp M\subseteq \Phi\}.$$
\end{thm}

Let us give an example of a set of prime ideals which is not
coherent. This is based on an example from \cite{BIK} and provides a
counterexample to Hovey's classification of thick subcategories closed
under direct sums in \cite[Theorem~5.2]{H}.

\begin{exm}\label{ex}
Let $k$ be a field and $A=k[\![x,y]\!]$. Then $\Phi=\Spec
A\setminus\{\frm\}$ with $\frm=(x,y)$ is not coherent. To see this,
let $$0\to A\to E(A)\to \bigoplus_{\height\frp=1} E(A/\frp)\to
E(A/\frm)\to 0$$ be a minimal injective resolution of $A$. Then the
morphism $$E(A)\to \bigoplus_{\height\frp=1} E(A/\frp)$$ cannot be
completed to an exact sequence lying in $\Inj_\Phi A$, because its
cokernel $E(A/\frm)$ does not belong to $\Inj_\Phi A$.
\end{exm}

It should be clear that one can construct such examples more generally
for commutative noetherian rings of Krull dimension at least two. On
the other hand, if $A$ is a Dedekind domain, then all subsets of
$\Spec A$ are coherent, because every factor module of an injective
$A$-module is injective and therefore the cokernel of a morphism
$I^0\to I^1$ between injective $A$-modules is up to isomorphism a
direct summand of $I^1$. We refer to Section~\ref{se:coh} for details
about coherent subsets.

The proof of Theorem~\ref{th:thick} uses an alternative description of
the support of a module. This involves the derived category $\bfD(\Mod
A)$ of $\Mod A$. Given two complexes $X$ and $Y$ of $A$-modules, we
write $X\lotimes_AY$ for their tensor product in $\bfD(\Mod A)$.
Note that for every $A$-module $M$, we have
$$\Supp M=\{\frp\in\Spec A\mid M\lotimes_A k(\frp)\neq 0\}.$$ The
following lemma is due to Foxby. The proof given here is inspired by
Neeman's work \cite[\S 2]{N}; see Proposition~\ref{pr:supp} for a more
general statement.
 
\begin{lem}[{\cite[Remark~2.9]{F}}]\label{le:supp}
Let $M$ be an $A$-module. Given a minimal injective resolution $I^*$
of $M$, we have $$\Supp M=\bigcup_{i\geq 0}\Ass I^i.$$
\end{lem}
\begin{proof}
First observe that we can pass from $M$ to the complex $I$, because
the morphism $M\to I$ induces an isomorphism in $\bfD(\Mod A)$. Fix a
prime ideal $\frp$.  Recall that each injective $A$-module $J$ admits
a unique decomposition
$$J=\bigoplus_{\frq\text{ prime}}\Gamma_\frq J$$ such that $\Ass
\Gamma_\frq J\subseteq\{\frq\}$ for all $\frq$.  We denote by
$\Gamma_\frp I$ the complex which is obtained from $I$ by taking in
each degree the component with associated prime $\frp$. To be precise,
$\Gamma_\frp I$ is the subcomplex of $I_\frp=I\otimes_AA_\frp$
supported at the closed point $\frp$. Note that the sequence $I\to
I_\frp\leftarrow \Gamma_\frp I$ of canonical morphisms is degreewise a
split epimorphism, followed by a split monomorphism.  In particular,
it induces an isomorphism
$$I\lotimes_Ak(\frp)\cong I_\frp\lotimes_Ak(\frp)\cong \Gamma_\frp
I\lotimes_Ak(\frp),$$ since $I'\lotimes_Ak(\frp)=0$ for the kernel
$I'$ of $I\to I_\frp$ and $I''\lotimes_Ak(\frp)=0$ for the cokernel
$I''$ of $\Gamma_\frp I\to I_\frp$.

Suppose first that $I\lotimes_Ak(\frp)\neq 0$. Then $\Gamma_\frp I\neq
0$ and therefore $\frp\in\Ass I^i$ for some $i$. Now suppose that
$I\lotimes_Ak(\frp)= 0$. Then $\Gamma_\frp I= 0$ by
\cite[Lemma~2.14]{N}. We want to conclude that
$\Gamma_\frp(I^i)=(\Gamma_\frp I)^i=0$ for all $i$. Here we need to
use the minimality of $I$. Recall that a complex $J$ of injective
$A$-modules is \emph{minimal} if for all $i$ the kernel of the
differential $J^i\to J^{i+1}$ is an essential submodule of $J^i$.  If
$J$ is minimal and $J^i=0$ for $i\ll 0$, then $H^iJ=0$ for all $i$
implies $J^i=0$ for all $i$. Observe that
$\Gamma_\frp$ preserves minimality. Thus $\Gamma_\frp I=0$ in
$\bfD(\Mod A)$ implies $\frp\not\in\Ass I^i$ for all $i$, because $I$
is minimal.
\end{proof}

\begin{lem}\label{le:coh}
Let $\Phi$ be a coherent subset of $\Spec A$.  Then 
$$\Supp^{-1}\Phi=\{M\in\Mod A\mid\Supp M\subseteq \Phi\}$$ is a thick
subcategory of $\Mod A$.
\end{lem}
\begin{proof} 
We consider the full subcategory $\C$ consisting of all $A$-modules $M$ which
fit into an exact sequence
$$0\to M\to I^0\to I^1\quad\text{with}\quad \Ass I^i\subseteq \Phi
\quad (i=0,1).$$ Without any assumptions on $\Phi$, it is clear that
$\C$ is an additive subcategory of $\Mod A$ which is closed under
taking kernels.  An application of the horseshoe lemma shows that $\C$
is closed under forming extensions. Next observe that $\C$ is closed
under taking cokernels. Here we use that $\Phi$ is coherent. By
definition, the cokernel of a morphism between injective modules in
$\C$ belongs to $\C$. A standard argument then shows that this
property extends to arbitrary morphisms in $\C$.  It follows from
Lemma~\ref{le:supp} that $\Supp^{-1}\Phi=\C$, and therefore
$\Supp^{-1}\Phi$ is thick.
\end{proof}

\begin{lem}\label{le:thick}
Let $\C$ be a subcategory of $\Mod A$ which is thick and closed under
taking direct sums.  Then the injective envelope $E(M)$ belongs to
$\C$ for every $M$ in $\C$.
\end{lem}
\begin{proof}
Fix $M$ in $\C$. First observe that $\Tor_i^A(M,N)$ belongs to $\C$
for every $A$-module $N$ and every integer $i$.  This is clear,
because for any projective resolution $P$ of $N$, the complex
$M\otimes_AP$ and therefore its cohomology lies in $\C$.  Given
$\frp\in\Supp M$, it follows that $k(\frp)$ belongs to $\C$, since
$\Tor_i^R(M,k(\frp))$ is a direct sum of copies of $k(\frp)$. Then
Lemma~\ref{le:inj} implies that $E(A/\frp)$ belongs to $\C$, and we
conclude from Lemma~\ref{le:supp} that $E(M)$ belongs to $\C$.
\end{proof}

\begin{proof}[Proof of Theorem~\ref{th:thick}]
Let $\Phi$ be a coherent subset of $\Spec A$.  Then the subcategory
$\Supp^{-1}\Phi$ is thick and closed under taking direct sums, by
Lemma~\ref{le:coh}.  Clearly, we have
$$\Supp(\Supp^{-1}\Phi)=\Phi.$$ Now let $\C$ be a subcategory of $\Mod A$,
which is thick and closed under taking direct sums. Let
$\Phi=\Supp\C$. First observe that $\Inj_\Phi A\subseteq\C$, by
Lemma~\ref{le:thick}. We claim that $\Phi$ is coherent. In deed, each
morphism $I^0\to I^1$ in $\Inj_\Phi A$ can be completed to an exact sequence
$I^0\to I^1\to I^2$ in $\Inj_\Phi A$ by taking for $I^2$ the injective
envelope of a cokernel of $I^0\to I^1$. Next we claim that
$$\Supp^{-1}(\Supp\C)=\C.$$ The inclusion
$\Supp^{-1}(\Supp\C)\supseteq\C$ is clear. Now suppose that $M$ is a
module contained in $\Supp^{-1}(\Supp\C)$ and choose a minimal
injective resolution $I^*$.  Then $\Supp I^i\subseteq\Supp\C$ for all
$i$. Thus $I^0$ and $I^1$ belong to $\C$, and we conclude that $M$
belongs to $\C$.
\end{proof}

The classification of thick subcategories specializes to Gabriel's
classification of localizing subcategories.

\begin{cor}[{\cite[p.~ 425]{G}}]
The map $\C\mapsto\Supp\C$ induces a bijection between the set
of localizing subcategories of $\Mod A$ and the set of specialization
closed subsets of $\Spec A$. 
\end{cor}

\section{Coherent subsets of $\Spec A$}\label{se:coh}

In this section we collect some basic properties of coherent subsets
of $\Spec A$.  Let us fix some notation. Given a multiplicatively
closed subset $S$ of $A$, let $\pi\colon A\to S^{-1}A$ denote the
localization.  Then we identify $\Spec S^{-1}A$ via $\pi^{-1}$ with
the subset of all prime ideals $\frp$ of $A$ satisfying
$S\cap\frp=\emptyset$.

\begin{prop}\label{pr:coh}
Let $\Phi$ be a subset of $\Spec A$.
\begin{enumerate}
\item Let $(\Phi_i)$ be a family of coherent subsets of $\Spec A$. Then
$\bigcap_i \Phi_i$ is coherent.
\item If $\Phi$ is specialization closed, then $\Phi$ is coherent.
\item If $\frq\subseteq\bigcup_{\frp\in\Phi}\frp$ implies
$\frq\in\Phi$ for every prime ideal $\frq$, then $\Phi$ is coherent.
\item The subset $\Phi$ is coherent if and only if $\Phi\cap\Spec A_\frp$ is
a coherent subset of $\Spec A_\frp$ for each prime ideal $\frp$.
\end{enumerate}
\end{prop}
\begin{proof}
We use that $\Phi$ is coherent if and only if the cokernel $C$ of each
morphism $I^0\to I^1$ between injective $A$-modules with $\Ass
I^i\subseteq \Phi$ ($i=0,1$) satisfies $\Ass C\subseteq \Phi$.

(1) Clear.

(2) The assumption on $\Phi$ implies that for each pair $N\subseteq M$
of $A$-modules with $\Ass M\subseteq \Phi$, we have that $\Ass
M/N\subseteq \Phi$.

(3) The set
$$S=A\setminus\bigcup_{\frp\in\Phi}\frp=\bigcap_{\frp\in
\Phi}A\setminus\frp$$ is multiplicatively closed. The assumption on
$\Phi$ implies that the localization $A\to S^{-1}A$ identifies all
injective $S^{-1}A$-modules with the injective $A$-modules $I$
satisfying $\Ass I\subseteq\Phi$.

(4) We write $\Phi_\frp=\Phi\cap\Spec A_\frp$ for each prime ideal
$\frp$. Suppose first that $\Phi$ is coherent and fix a prime ideal
$\frp$. Let $I^0\to I^1$ be a map in $\Inj_{\Phi_\frp} A$. There exist
an exact sequence $I^0\to I^1\to I^2$ in $\Inj_\Phi A$ and
localization at $\frp$ induces an exact sequence in $\Inj_{\Phi_\frp}
A$. Thus $\Phi_\frp$ is coherent. Now suppose that $\Phi$ is not
coherent. It follows that there exists an exact sequence $$I^0\to
I^1\to C\to 0$$ of $A$-modules with $I^i\in \Inj_\Phi A$ ($i=0,1$) but
$\Ass C\not\subseteq \Phi$. Let $\frp\in\Ass C\setminus \Phi$.  Then
we localize at $\frp$ and obtain an exact sequence $$I^0_\frp\to
I^1_\frp\to C_\frp\to 0$$ of $A_\frp$-modules with $I^i_\frp\in
\Inj_{\Phi_\frp} A$ ($i=0,1$) but
$$\Ass C_\frp=(\Ass C)\cap\Spec A_\frp\not\subseteq \Phi_\frp.$$ Thus
$\Phi_\frp$ is a subset of $\Spec A_\frp$ which is not coherent.
\end{proof}

\begin{rem}
(1) A subset $\Phi$ of $\Spec A$ satisfies the condition (3)
of Proposition~\ref{pr:coh} if and only if it is of the form $\Spec
S^{-1}A$ for some multiplicatively closed subset $S$. 

(2) The union of two coherent subsets need not to be coherent. For
instance, Example~\ref{ex} provides a subset which is not coherent but
Zariski open. Each Zariski open subset $U$ can be written as the
finite union of basic open subsets. However, a basic open set is
coherent because it is of the form $\Spec S^{-1}A$ for some
multiplicatively closed subset $S$.
\end{rem}

\begin{cor}\label{co:coh}
If the Krull dimension of $A$ is at most one, then every subset of
$\Spec A$ is coherent. 
\end{cor}
The converse of this statement is proved in the appendix of this paper.
\begin{proof}
We may assume that $A$ is local, by part (4) of
Proposition~\ref{pr:coh}, and we denote by $\frm$ the maximal ideal.
Let $\Phi$ be a subset of $\Spec A$. If $\Phi$ contains $\frm$, then
$\Phi$ is specialization closed and therefore coherent, by part (2) of
Proposition~\ref{pr:coh}. If $\frm$ is not contained in $\Phi$, then
all prime ideals in $\Phi$ are minimal and therefore the prime
avoidance theorem implies that the condition in part (3) of
Proposition~\ref{pr:coh} is satisfied. Thus $\Phi$ is coherent.
\end{proof}

Given a prime ideal $\frp$ of $A$, let
$$V(\frp)=\{\frq\in\Spec A\mid \frp\subseteq \frq\}\quad
\text{and}\quad \varLambda(\frp)=\{\frq\in\Spec A\mid \frq\subseteq
\frp\}.$$ Subsets of the from $V(\frp)$ and $\varLambda(\frp)$ are
coherent. They can be used to build new coherent subsets.

\begin{cor}
Let $\Phi$ and $\Psi$ be subsets of $\Spec A$ and suppose that $\Psi$ is
finite.  Then
$$\bigcup_{\frp\in\Phi,\frq\in\Psi}
V(\frp)\cap \varLambda(\frq)$$ is a coherent subset of $\Spec A$.
\end{cor}
\begin{proof}
We can express the set as the intersection of two coherent subsets:
\begin{equation*}
\bigcup_{\frp\in\Phi,\frq\in\Psi} V(\frp)\cap
\varLambda(\frq)=\big(\bigcup_{\frp\in\Phi} V(\frp)\big)\cap
\big(\bigcup_{\frq\in\Psi} \varLambda(\frq)\big).\qedhere
\end{equation*}
\end{proof}

\begin{exm}
Let $\frp_1,\frp_2$ be prime ideals. If $\frp_1\subseteq\frp_2$, then
$\{\frq\mid\frp_1\subseteq\frq\subseteq\frp_2\}$ is coherent. If
$\frp_1\not\subseteq\frp_2$ and $\frp_2\not\subseteq\frp_1$, then
$\{\frp_1,\frp_2\}$ is coherent. In both cases the set is of the form
$$\bigcup_{\frp,\frq\in\Phi} V(\frp)\cap
\varLambda(\frq)\quad\text{with}\quad\Phi=\{\frp_1,\frp_2\}.$$
\end{exm}

\section{Support of complexes}

Let $X$ be a complex of $A$-modules. Following \cite{F}, the
\emph{support} of $X$ is by definition
$$\Supp X=\{\frp\in\Spec A\mid X\lotimes_A k(\frp)\neq 0\}.$$ We use
an alternative description of the support of $X$.  A complex $I$ of
injective $A$-modules together with a quasi-isomorphism $X\to I$ is
called a \emph{minimal K-injective resolution} of $X$, if $I$ is
\emph{K-injective} (that is, every morphism from an acyclic complex to
$I$ is null homotopic) and $I$ is \emph{minimal} (that is, for all $i$
the kernel of the differential $I^i\to I^{i+1}$ is an essential
submodule of $I^i$).

One can show that each complex of $A$-modules admits a minimal
K-injective resolution; see \cite[Theorem~4.5]{S} or
\cite[Application~2.4]{BN} for the existence of a K-injective
resolution and \cite[Proposition~B.2]{K} for the minimality.  Note
that each acyclic and K-injective complex is null homotopic. Moreover,
if a minimal complex $I$ of injective $A$-modules is null homotopic,
then $I^i=0$ for all $i$.

The next proposition is the obvious generalization of
Lemma~\ref{le:supp} from modules to complexes of modules. The proof
requires only minor modifications; it follows closely \cite[\S 2]{N}.

\begin{prop}\label{pr:supp}
Let $X$ be a complex of $A$-modules and $X\to I$ a minimal K-injective
resolution of $X$.  Then we have $$\Supp X=\bigcup_{i\in\bbZ}\Ass
I^i.$$
\end{prop}
\begin{proof}
First observe that we can pass from $X$ to the complex $I$, because
the morphism $X\to I$ induces an isomorphism in $\bfD(\Mod A)$. Fix a
prime ideal $\frp$.  Recall that each injective $A$-module $J$ admits
a unique decomposition
$$J=\bigoplus_{\frq\text{ prime}}\Gamma_\frq J$$ such that $\Ass
\Gamma_\frq J\subseteq\{\frq\}$ for all $\frq$.  We denote by
$\Gamma_\frp I$ the complex which is obtained from $I$ by taking in
each degree the component with associated prime $\frp$. To be precise,
$\Gamma_\frp I$ is the subcomplex of $I_\frp=I\otimes_AA_\frp$
supported at the closed point $\frp$. Note that the sequence $I\to
I_\frp\leftarrow \Gamma_\frp I$ of canonical morphisms is degreewise a
split epimorphism, followed by a split monomorphism.  In particular,
it induces an isomorphism
$$I\lotimes_Ak(\frp)\cong I_\frp\lotimes_Ak(\frp)\cong \Gamma_\frp
I\lotimes_Ak(\frp),$$ since $I'\lotimes_Ak(\frp)=0$ for the kernel
$I'$ of $I\to I_\frp$ and $I''\lotimes_Ak(\frp)=0$ for the cokernel
$I''$ of $\Gamma_\frp I\to I_\frp$.

Suppose first that $I\lotimes_Ak(\frp)\neq 0$. Then $\Gamma_\frp I\neq
0$ and therefore $\frp\in\Ass I^i$ for some $i$. Now suppose that
$I\lotimes_Ak(\frp)= 0$. Then $\Gamma_\frp I= 0$ by
\cite[Lemma~2.14]{N}, that is, $\Gamma_\frp I$ is acyclic. We want to
conclude that $\Gamma_\frp(I^i)=(\Gamma_\frp I)^i=0$ for all $i$. Here
we need to use the minimality of $I$.  We observe that $\Gamma_\frp$
preserves minimality. Also, $\Gamma_\frp I$ is a K-injective complex
of injective $A$-modules.  Thus $\Gamma_\frp I=0$ in $\bfD(\Mod A)$
implies $\frp\not\in\Ass I^i$ for all $i$, because $I$ is minimal.
\end{proof}

%Recall from \cite{N} that the map sending a subcategory $\C$ of
%$\bfD(\Mod A)$ to $$\Supp\C=\bigcup_{X\in\C}\Supp X$$ induces a
%bijection between the set of localizing subcategories of $\bfD(\Mod
%A)$ and the set of subsets of $\Spec A$.

\begin{thm}
For a subset $\Phi$ of $\Spec A$ the following conditions are equivalent:
\begin{enumerate}
\item $\Phi$ is coherent.
\item For every complex $X$ of $A$-modules we have
$$\Supp X\subseteq \Phi \quad\Longleftrightarrow\quad \Supp H^iX\subseteq \Phi\text{
for all }i\in\bbZ.$$
\item For every complex $X$ of $A$-modules we have
$$\Supp X\subseteq \Phi \quad\Longrightarrow\quad \Supp H^iX\subseteq \Phi\text{
for all }i\in\bbZ.$$
\end{enumerate}
\end{thm}
\begin{proof}
(1) $\Rightarrow$ (2): Suppose that $\Phi$ is coherent.  We use that the $A$-modules
$M$ with $\Supp M\subseteq \Phi$ form a thick subcategory, by
Lemma~\ref{le:coh}. Now fix a complex $X$ of $A$-modules. If $\Supp
X\subseteq \Phi$, then we have a quasi-isomorphic complex $I$ of
injective $A$-modules with $\Ass I^i\subseteq \Phi$ for all $i$, by
Proposition~\ref{pr:supp}. It follows that
$$\Supp H^iX=\Supp H^iI\subseteq \Phi$$ for all $i$. Now let $\Supp
H^*X\subseteq \Phi$ and $\frp\in\Supp X$.  A d\'evissage argument
shows that $\Supp H^*(X\lotimes_A Y)\subseteq \Phi$ for every complex
$Y$. In particular,
$$\{\frp\}=\Supp H^*(X\lotimes_Ak(\frp))\subseteq \Phi$$ because
$X\lotimes_Ak(\frp)$ is a direct sum of shifted copies of $k(\frp)$. Thus
$\Supp X\subseteq \Phi$.

(2) $\Rightarrow$ (3): Clear.

(3) $\Rightarrow$ (1): Let $I^0\to I^1$ be a morphism of injective
$A$-modules with $\Ass I^i\subseteq \Phi$ for $i=0,1$. Viewing $I$ as
a complex, we have $\Supp I\subseteq \Phi$ by
Proposition~\ref{pr:supp}, and therefore $\Supp H^0I\subseteq
\Phi$. It follows from Lemma~\ref{le:supp} that we can complete
$I^0\to I^1$ to an injective resolution $$0\to H^0I\to I^0\to I^1\to
I^2\to \cdots$$ of $H^0I$ with $\Ass I^i\subseteq \Phi$ for all
$i$. Thus $\Phi$ is coherent.
\end{proof}

\subsection*{Acknowledgements} 
A pleasant collaboration with Dave Benson and Srikanth Iyengar has
been the starting point of this paper; in particular the
Example~\ref{ex} is taken from \cite{BIK}.  I would like to thank
Amnon Neeman and Torsten Wedhorn for a number of helpful comments on
this work.

\begin{appendix}
\section{Noncoherent subsets of $\Spec A$}
\begin{center}
{\sc By Srikanth Iyengar}
\end{center}

\medskip
In this appendix, we establish the converse of Corollary~\ref{co:coh}.  To this end, we
recall some standard notions from commutative algebra; this serves also to fix notation.

\begin{defn}
Let $A$ be a commutative noetherian ring, $\fra$ an ideal in $A$, and let $M$ be an $A$-module.
The \emph{$\fra$-depth} of $M$ is the number
\[
\depth_A(\fra,M) = \inf\{n\mid \Ext^n_A(A/\fra,M)\ne 0\}\,.
\]  
This invariant of $M$ can also be detected from its Koszul homology on a finite set of
elements generating $\fra$, and also its local cohomology with respect to $\fra$; see, for
instance, \cite[Theorem 2.1]{FI}.  When $M$ is finitely generated and $\fra M\ne M$, this number
coincides with the length of the longest $M$-regular sequence in $\fra$; see \cite[Theorem 28]{Ma}.

As usual, if $A$ is local, with maximal ideal $\frm$, we write $\depth_AM$ for the
$\frm$-depth of $M$, and call it the \emph{depth} of $M$.
\end{defn}

We record the following standard properties of depth, for ease of reference.

\begin{lem}
\label{depth}
Let $\fra$ be an ideal in a commutative noetherian ring $A$, and let $M$ be an $A$-module. 
The following statements hold.
\begin{enumerate}
\item One has an equality, $\depth_A(\sqrt{\fra},M)=\depth_A(\fra,M)$.
\item With $I^*$ the minimal injective resolution of $M$, for each prime ideal $\frp$ one has 
\[
\depth_{A_\frp}M_\frp = \inf\{n\mid \text{$E(A/\frp)$ is a direct summand in $I^n$}\}\,.
\]
\item If $A\to B$ is a homomorphism of rings and $N$ is a $B$-module, then viewing $N$ as
  an $A$-module by restriction of scalars, one has
\[
\depth_A(\fra,N)=\depth_B(\fra B,N)\,.
\]
\end{enumerate}
\end{lem}

\begin{proof}
For (1) see, for instance, \cite[Proposition~2.11]{FI}, while (3) is
evident, if one computes depth using Koszul complexes, or via local
cohomology; see \cite[Theorem~2.1]{FI}.  Part (2) holds as
$(I^*)_\frp$ is the minimal injective resolution of $M_\frp$ over
$A_\frp$, so the complex $\Hom_{A_\frp}(A_\frp/\frp
A_\frp,(I^*)_\frp)$, whose first nonzero cohomology module occurs in
degree $\depth_{A_\frp}M_\frp$, has zero differential.
\end{proof}

We need the following result of Auslander and Buchsbaum, which is implicit in \cite{AB}.

\begin{prop}
\label{pr:depth2}
Let $A$ be a commutative noetherian ring with $\dim A\ge 1$. There
exists a prime $\frp$ in $\Spec A$ such that 
$$\depth A_\frp = \dim A_\frp = \dim A-1.$$
\end{prop}

\begin{proof}
  The proof uses an induction on $\dim A$. When $\dim A=1$, for $\frp$ one may take any
  minimal prime of $A$.  This is the basis of the induction.

  Assume $\dim A\ge 2$, and let $\frm$ be the maximal ideal of $A$. Since $\Ass A$ is
  finite, the prime avoidance theorem implies that the following set in nonempty:
\[
\frm \setminus \bigcup_{\overset{\frp\in \Ass A}{\frp\ne \frm}}\frp\,.
\]
Choose an element $x$ in it. One then has that $\dim(A/Ax) =\dim A-1$, so the induction
hypothesis yields a prime $\frp$ in $A$, containing $x$, such that
\[
\depth (A_\frp/A_\frp x) =  \dim (A_\frp/A_\frp x) = \dim (A/Ax) -1 = \dim A -2\,.
\]
Observe that $\frp\ne \frm$, so the choice of $x$ ensures that it is a nonzero divisor in
$A_\frp$. Thus one has $\depth A_\frp = \dim A_\frp = \dim A - 1$. This completes the
induction argument.
\end{proof}

The gist of the result below is well-known; we provide a proof for
lack of a suitable reference for this formulation.

\begin{thm}
\label{thm:depth2}
Let $A$ be a commutative noetherian ring with $\dim A\ge 2$. There exists an $A$-module
$M$ and a prime $\frp$ in $\Spec A$ such that $\depth_{A_\frp}M_\frp =\max\{2,\dim
A-1\}$.
\end{thm}

\begin{proof}
  When $\dim A\ge 3$, we apply Proposition~\ref{pr:depth2}, and take for $\frp$ any prime
  such that $\depth A_\frp = \dim A-1$ and set $M=A$.  Henceforth, we assume $\dim A=2$.
  Choosing a prime ideal $\frp$ in $A$ with $\dim A_\frp=2$, and replacing $A$ with
  $A_\frp$, we may assume also that $A$ is local; the goal then is to find an $A$-module
  $M$ such that $\depth_AM=2$.  At this point, one may refer to, for instance, Hochster's
  article~\cite{Ho}, especially Section~3. We provide details, for completeness.

  Let $\widehat A$ be the completion of $A$ at its maximal ideal, say
  $\frm$.  One has $\dim \widehat A=2$, so there exists a prime ideal
  $\fra$ in $\widehat A$ with $\dim (\widehat A/\fra)=2$. Consider the
  canonical homomorphisms of rings $A\to \widehat A\to \widehat
  A/\fra$. Observer that $\frm(\widehat A/\fra)$ is the maximal ideal
  of $\widehat A/\fra$, so, by Lemma~\ref{depth}(3), for any module
  $M$ over $\widehat A/\fra$, one has
\[
\depth_AM=\depth_{\widehat A/\fra}M\,.
\]
Replacing $A$ with $\widehat A/\fra$ one may assume $A$ is a complete local domain, with
$\dim A=2$.

Let $B$ be the integral closure of $A$ in its field of fractions; the
conditions on $A$ imply that $B$ is finite as an $A$-module, by
\cite[Corollary~2, pp. 234]{Ma}, so also a two dimensional noetherian
ring, by \cite[Theorem~20]{Ma}. The finiteness of the extension
$A\subseteq B$ implies that there exists a prime ideal $\frq$ in $B$
with $\dim B_\frq=2$ and $\frq\cap A=\frm$, the maximal ideal of $A$;
see \cite[Theorem~5(iii)]{Ma}. The choice of $\frq$ ensures that
$\sqrt{\frm B_\frq}=\frq B_\frp$, so part (3) and (2) of
Lemma~\ref{depth} yield the first and second equalities below:
  \begin{align*}
    \depth_A B_\frq     &= \depth_{B_\frq}(\frm B_\frq,B_\frq)\\
                        &= \depth B_\frq \\
                        &= \dim B_\frq \\
                        &=2
  \end{align*}
  The penultimate equality holds because $B_\frq$ is Cohen-Macaulay, by Serre's criterion
  for normality; see \cite[Theorem~39]{Ma}. The $A$-module $B_\frq$ has thus the desired
  depth.
\end{proof}

The result below is a perfect converse to Corollary~\ref{co:coh}.

\begin{cor}
  If $A$ is a commutative noetherian ring with $\dim A\ge 2$, then there exists a subset
  $\Phi$ of $\Spec A$ that is not coherent.
\end{cor}

\begin{proof}
  By Theorem~\ref{thm:depth2}, there exists an $A$-module $M$ and a prime ideal $\frp$
  in $A$ such that $d=\depth_{A_\frp}M_\frp\ge 2$. Let $I^*$ be a minimal injective
  resolution of $M$, and set
\[
\Phi = \Ass_A(I^{d-2}) \cup \Ass_A (I^{d-1})
\]
Observe that $\frp$ is in $\Ass_A(I^d)$ but not in $\Phi$, by
Lemma~\ref{depth}(2). Thus, the set $\Phi$ is not coherent, as the
associated primes of $\Coker(I^{d-2}\to I^{d-1})$ coincide with those
of $I^d$.
\end{proof}
\end{appendix}


\begin{thebibliography}{99}
%
\bibitem{AB} {\sc M. Auslander and D. Buchsbaum:} Homological dimension in noetherian
  rings. II.  Trans. Amer. Math. Soc. \textbf{88} (1958), 194--206.

%
\bibitem{BIK}{\sc D.\ Benson, S.\ Iyengar, and H.\ Krause:} Local
cohomology and support for triangulated categories. arXiv:math.KT/0702610. 
%
\bibitem{BN}{\sc M. B\"okstedt and A. Neeman:} Homotopy limits in
triangulated categories. Compositio Math. {\bf 86} (1993) 209--234.
%
\bibitem{B0}{\sc N. Bourbaki:} Alg\`ebre commutative. Hermann, Paris, 1968. 
%
\bibitem{B}{\sc N. Bourbaki:} Alg\`ebre commutative. Chapitre 10. Masson, Paris, 1998. 
%
\bibitem{F}{\sc H.-B. Foxby:} Bounded complexes of flat modules.
J. Pure Appl. Algebra \textbf{15} (1979), 149--172.
%
\bibitem{FI}{\sc H.-B. Foxby and S. Iyengar:} Depth and amplitude for unbounded complexes.
  Commutative algebra (Grenoble/Lyon, 2001), 119--137, Contemp. Math., {\bf 331},
  Amer. Math. Soc., Providence, RI, 2003.
%
\bibitem{G}{\sc P. Gabriel:} Des cat\'egories
ab\'eliennes. Bull. Soc. Math. France {\bf 90} (1962), 323--448.
%
\bibitem{GP}{\sc G. Garkusha and M. Prest:} Classifying Serre
subcategories of finitely presented modules. Preprint (2006).
%
\bibitem{Ho}{\sc M. Hochster:} Cohen-Macaulay modules. Conference on Commutative Algebra
  (Lawrence, Kansas 1972), 120--152, Lecture Notes Math. \textbf{311}, Springer-Verlag, 1973.
%
\bibitem{H}{\sc M. Hovey:} Classifying subcategories of
modules. Trans. Amer. Math. Soc. \textbf{353} (2001), 3181--3191.
%
\bibitem{K}{\sc H. Krause:} The stable derived category of a noetherian scheme.
  Compos. Math.  \textbf{141} (2005), 1128--1162.
%
\bibitem{Ma}{\sc H. Matsumura:} Commutative algebra. Second edition. Mathematics Lecture
  Note Series, \textbf{56}, Benjamin/Cummings Publishing Co., Inc., Reading, Mass., 1980.
%
\bibitem{N}{\sc A. Neeman:} The chromatic tower of $D(R)$. Topology
\textbf{31} (1992), 519--532.
%
\bibitem{S}{\sc N. Spaltenstein:} Resolutions of unbounded
complexes. Compositio Math. {\bf 65} (1988), 121--154.
%
\bibitem{T}{\sc R. Takahashi:} Classifying subcategories of modules
over commutative noetherian rings. Preprint (2006).
%
\end{thebibliography}
\end{document}